\documentclass[11pt,a4paper]{amsart}

\date{October 31, 2005}

\allowdisplaybreaks

\AtBeginDocument{\baselineskip=23pt} 

\newcommand{\htr}{H_{\alpha}} 
\newcommand{\ma}{M_{\alpha}}
\newcommand{\mua}{d\mu_\alpha}
\newcommand{\jna}{j_n^{\alpha}}
\newcommand{\jka}{j_k^{\alpha}}
\newcommand{\cna}{c_n^{\alpha}}
\newcommand{\cka}{c_k^{\alpha}}
\newcommand{\Sna}{S_n^{\alpha}}

\newcommand{\spa}{\operatornamewithlimits{span}} 

\newcommand{\calLpsa}{\mathcal{L}^p_{s,\alpha}} 

\newcommand{\lpnorm}[3]{\|{#1}\|_{L^{#2}({#3})}} 
\newcommand{\Biglpnorm}[3]{\Big\|{#1}\Big\|_{L^{#2}({#3})}} 
\newcommand{\calnorm}[1]{\|{#1}\|_{\calLpsa}} 

\newtheorem{theor}{Theorem}
\newtheorem{propo}{Proposition}
\newtheorem{lem}{Lemma}

\theoremstyle{remark}
\newtheorem*{rem}{Remark}

\author[Betancor]{Jorge J. Betancor}
\address[Betancor]{Departamento de An\'alisis Matem\'atico \\
         Universidad de La Laguna \\
         38271 La Laguna (Tenerife) \\
         Islas Canarias, Spain}
\email{jbetanco@ull.es}

\author[Ciaurri]{\'Oscar Ciaurri}
\address[Ciaurri, Varona]{Departamento de Matem\'aticas y Computaci\'on \\
         Universidad de La Rioja                       \\
         26004 Logro\~no, Spain}
\email{oscar.ciaurri@dmc.unirioja.es, jvarona@dmc.unirioja.es}
\urladdr{http://www.unirioja.es/dptos/dmc/jvarona/welcome.html}

\author[Mart\'{\i}nez]{Teresa Mart\'{\i}nez}
\address[Mart\'{\i}nez, Torrea]{Departamento de Matem\'aticas \\
         Facultad de Ciencias \\
         Universidad Aut\'o\-no\-ma de Madrid \\
         28049 Madrid, Spain}
\email{teresa.martinez@uam.es, joseluis.torrea@uam.es}
\urladdr{http://www.uam.es/personal\_pdi/ciencias/torrea}

\author[P\'erez]{Mario P\'erez}
\address[P\'erez]{Departamento de Matem\'aticas \\
         Universidad de Zara\-go\-za \\
         50009 Zara\-go\-za, Spain}
\email{mperez@unizar.es}
\urladdr{http://www.unizar.es/analisis\_matematico/mperez/mperez.html}

\author[Torrea]{Jos\'e Luis Torrea}

\author[Varona]{Juan L. Varona}

\thanks{Research supported by grants BFM2002-04013-C02-02 and
BFM2003-06335-C03-03 of the~DGI}



\begin{document}

\title{Heat and Poisson semigroups for Fourier-Neumann expansions}

\begin{abstract}
Given $\alpha > -1$, consider the second order differential
operator in $(0,\infty)$
$$
  L_\alpha f \equiv \left(x^2\, \frac{d^2}{dx^2}
    + (2\alpha+3)x\, \frac{d}{dx} + x^2 + (\alpha+1)^2 \right)(f),
$$
which appears in the theory of Bessel functions. The purpose of this paper
is to develop the corresponding harmonic analysis taking $L_\alpha$ as
the analogue to the classical Laplacian. Namely we study the boundedness
properties of the heat and Poisson semigroups. These boundedness properties
allow us to obtain some convergence results that can be used to solve the Cauchy
problem for the corresponding heat and Poisson equations.
\end{abstract}

\maketitle

Given $\alpha > -1$, we shall consider the second order operator on functions defined on $(0,\infty)$
$$
  L_\alpha f \equiv \left(x^2\, \frac{d^2}{dx^2}
    + (2\alpha+3)x\, \frac{d}{dx} + x^2 + (\alpha+1)^2 \right)(f).
$$
This operator appears in the theory of Bessel functions (see~\cite{Wat}).
It is selfadjoint with respect to the measure $\mua(x) = x^{2\alpha+1}\,dx$. It is well known that the functions
\begin{equation*}
  \jna(x) = \sqrt{2(\alpha+2n+1)}\,J_{\alpha+2n+1}(x) x^{-\alpha-1},
  \qquad n = 0,1,2,\dots\,
\end{equation*}
where $J_\nu$ stands for the Bessel function of the first kind of order~$\nu$, are eigenfunctions of the operator~$ L_\alpha$.
In fact 
$$
  L_\alpha \jna = (\alpha+2n+1)^2 \jna, \quad n=0,1,\dots
$$
see \cite[\S\,5.73, p.~158]{Wat}. By using some classical formulas for the Bessel functions $J_\nu$, it is easy to check that $\{\jna\}_{n=0}^\infty$ is an orthonormal system in
$L^2((0,\infty),\mua)$ ($L^2(\mua)$ from now on). 

Consider the so-called modified Hankel transform $H_\alpha$, that is 
\begin{equation}
\label{eq:htr}
  \htr f(x) = \int_0^\infty \frac{J_\alpha(xy)}{(xy)^\alpha} \, f(y)
    y^{2\alpha+1} \, dy,
  \quad x>0.
\end{equation}
Since it is known that $\htr \jna$ is supported on $[0,1]$, and $\htr$ is an isometry on $L^2(\mua)$, this system $\{\jna\}_{n=0}^\infty$ is not complete in $L^2(\mua)$.  On the other hand the subspace
$B_{2,\alpha} = \overline{\spa\{\jna\}_{n=0}^\infty}$ (closure in $L^2(\mua)$)
can be identified with the space
$\{ f \in L^2(\mua): \ma f = f\} = \ma(L^2(\mua))$,
where $\ma$ is the multiplier defined by $\htr(\ma f) = \chi_{[0,1]}\htr f$;
see \cite{V,CiStV}.

Along this paper we shall consider the operator $L_\alpha$ as a positive selfadjoint operator defined in the Hilbert space $B_{2,\alpha}$. Then, its heat and Poisson semigroups $W_t = e^{-tL_\alpha}$ and $P_t = e^{-t\sqrt{L_\alpha}}$ can be defined in a spectral way as
\begin{align*}
   W_t f &= \sum_{n=0}^\infty e^{-t(\alpha+2n+1)^2} c_n \jna,
   \\
   P_t f &= \sum_{n=0}^\infty e^{-t(\alpha+2n+1)} c_n \jna,
\end{align*}
for $f \in B_{2,\alpha}$ given by $f=\sum_{n=0}^\infty c_n \jna$.
Using the ideas in \cite{Stein} we could also define the Poisson semigroup $P_t$ by the following subordination formula:
\[
  P_tf(x)
  = \frac{1}{\sqrt{2\pi}} \int_0^\infty t e^{-t^2/4s}W_sf(x) s^{-3/2}\,ds,
\]
which can be derived from the well known identity
\[
   e^{-t\sqrt{\gamma}} =
   \frac{1}{\sqrt{2\pi}} \int_0^\infty t e^{-t^2/(4s)}e^{-s\gamma} s^{-3/2}\,ds.
\]
Analogously the formula
$s^{-\lambda} = \frac{1}{\Gamma(\lambda)}\int_0^\infty t^{\lambda-1} e^{-ts}\,dt$
suggests the definition of the Riesz potentials either  as
\begin{equation}
\label{eq:frac}
  L_\alpha^{-\lambda}f(x) =
  \frac{1}{\Gamma(\lambda)}\int_0^\infty t^{\lambda-1} W_tf(x)\,dt,
\end{equation}
or
\begin{equation}
\label{eq:Poisson}
   L_\alpha^{-\lambda/2}f(x)
   = \frac{1}{\Gamma(\lambda)} \int_0^\infty t^{\lambda-1} P_tf(x) \,dt.
\end{equation}

The operators considered above, $W_t$, $P_t$, and $L_\alpha^{-\lambda}$, are clearly bounded in $B_{2,\alpha}$. A natural question is to analyse the boundedness of these operators in the spaces $B_{p,\alpha}$, $1<p<\infty$, defined as the closure in $L^p(\mua)$ of the space $\spa\{\jna\}_{n=0}^\infty$.
The first requirement is that $\jna \in L^p(\mua)$ for every~$n$. 
By using well known estimates for the Bessel functions (see~\eqref{eq:zero} and~\eqref{eq:infty}), this implies $p > \max\{1,\frac{4(\alpha+1)}{2\alpha+3}\}$. Moreover, if the Fourier coefficients $\cna(f)$
must exist for every $f \in L^p(\mua)$, this requires $\jna \in L^{p'}(\mua)$ for every~$n$, where $1/p+1/p'=1$. This is equivalent to $p < \frac{4(\alpha+1)}{2\alpha+1}$ if $\alpha\ge-1/2$, and $p < \infty$, when $-1 < \alpha < -1/2$. This observation leads us to restrict our study to the space $B_{p,\alpha}$ with $p \in (p_0(\alpha), p_1(\alpha))$ where 
$$
  p_0(\alpha) =
  \begin{cases}
    \frac{4(\alpha+1)}{2\alpha+3} & \text{if } \alpha \ge -1/2,\\
    1 & \text{if } -1 < \alpha<-1/2,
  \end{cases}
$$
and
$$
  p_1(\alpha) =
  \begin{cases}
    \frac{4(\alpha+1)}{2\alpha+1} & \text{if } \alpha \ge -1/2,\\
    \infty & \text{if } -1 < \alpha<-1/2.
  \end{cases}
$$
As in many other cases in the literature, technical reasons make convenient the change of parameter $r=e^{-t}$. With a small abuse of notation, we still use the notation $P_r$ and $W_r$ for the corresponding Poisson and heat semigroups. We prove, see Theorems~\ref{th:Poisson} and~\ref{th:heat}, that the operators $P_r$ and $W_r$ 
are uniformly bounded in $L^p(\mua)$, for $p \in (p_0(\alpha), p_1(\alpha))$.
As usual, the uniform boundedness produce the corresponding mean
convergence results, see Theorem~\ref{th:convLpVr}. 
Surprisingly, in this case the mean convergence allows us to prove the almost everywhere convergence, see Theorem~\ref{th:convaeVr}. This is due to the decay of the involved kernels.
We should mention that in the case $\alpha \ge -1/2$ the space $B_{p,\alpha}$ was characterized as $\ma(L^p(\mua))$, then some special results can be derived in this situation, see Theorems~\ref{th:convMaf} and~\ref{th:convaeMaf}. As a byproduct of the proof of Theorem~\ref{th:heat}, we find for every $(t,x)\in (0,\infty)\times(0,\infty)$ an expression for the function $W_t f(x)$, when $f\in B_{p,\alpha}$. In Theorem~\ref{th:DE} we prove that this function is infinitely differentiable with respect to both variables $t$ and $x$, and that it satisfies the heat equation $\left( \frac{\partial}{\partial t} + L_\alpha \right) W_t f(x) = 0$. The already mentioned convergence results give some solution for the corresponding Cauchy problem. Some applications to fractional integrals and potential spaces are also considered in the last section of the paper.

One could say that the spirit of this paper is nothing but developing a harmonic analysis associated to a second order differential operator in a parallel way to the classical Laplacian. This idea has the names of Muckenhoput~\cite{Mu} and Stein~\cite{Stein-2} as pioneer authors. In the last decade there was a big flourishing in this area and a relatively large number of papers appeared, see~\cite{MST, BMTU}.
 
\section{Technical results}

The Bessel functions satisfy the asymptotic formulas
(see, for instance,~\cite[Ch.~III, 3.1 (8), p.~40]{Wat} and~\cite[Ch.~VII, 7.21 (1),
p.~199]{Wat})
\begin{equation}
\label{eq:zero}
  J_\nu(x) = \frac{x^\nu}{2^\nu \Gamma(\nu+1)} + O(x^{\nu+2}),
    \qquad x \to 0^+,
\end{equation}
\begin{equation}
\label{eq:infty}
  J_\nu(x) = \left(\frac{2}{\pi x}\right)^{1/2}\left[ \cos\left(x-\frac{\nu\pi}2
    - \frac\pi4 \right) + O(x^{-1}) \right], \qquad x \to \infty.
\end{equation}

We shall also use the following estimates that can be found in~\cite{CiGPV,V}:
\begin{equation}
\label{eq:jnu}
   |J_\nu(x)| \le C x^{-1/4} \left( |x - \nu| + \nu^{1/3} \right)^{-1/4},
   \quad x \in (0,\infty),
\end{equation}
where $C$ is a positive constant independent of~$\nu$. 

\begin{lem}
\label{lem:normajn}
Let $\alpha > -1$ and $p_0(\alpha)<p<\infty$.
Then, $\{\jna\}_{n=0}^\infty \subseteq L^p(\mua)$ and
\begin{equation*}
  \lpnorm{\jna}{p}{\mua} \le C
    \begin{cases}
      n^{-(\alpha+1) + 2(\alpha+1)/p},
      &\text{ if } p < 4, \\[2pt]
      n^{-(\alpha+1)/2} (\log n)^{1/4},
      &\text{ if } p = 4, \\[2pt]
      n^{-(5/6+\alpha) + (6\alpha + 4)/(3p)},
      &\text{ if } p > 4.
    \end{cases}
\end{equation*}
\end{lem}

\begin{proof}
The assertion $\jna \in L^p(\mua)$ for every $n=0,1,2,\dots$
follows from~\eqref{eq:zero} and~\eqref{eq:infty}. Then,
estimates \eqref{eq:jnu} above show that $\lpnorm{\jna}{p}{\mua}$ is bounded above by a constant times the right hand side. For a similar expression, see~\cite{St-Guo}.
\end{proof}

\begin{lem}
\label{lem:convae}
Let $\alpha > -1$ and $p$ with $p_0(\alpha)<p< p_1(\alpha)$.
Then, for any $f \in L^p(\mua)$ the Fourier series $\sum_{n=0}^\infty \cna(f)\jna(x)$ converges absolutely for every $x \in (0,\infty)$. (Note that we do not assert that this convergence is to $f(x)$, not even almost everywhere.)
\end{lem}

\begin{proof}
Recall that
\begin{equation}
\label{eq:fourier}
  \cna(f) = \int_0^\infty f(y) \jna(y) y^{2\alpha+1} \, dy.
\end{equation}
It follows from Lemma~\ref{lem:normajn} that
$\jna \in L^{p'}(\mua)$; moreover,
$\lpnorm{\jna}{p'}{\mua} \le C n^\delta$
for some constant $\delta = \delta(p,\alpha)$. Thus,
by H\"older's inequality,
\begin{equation*}
   |\cna(f)|
   \le \lpnorm{f}p{\mua}
       \lpnorm{\jna}{p'}{\mua}
   \le C \lpnorm{f}p{\mua}\,n^\delta.
\end{equation*}
Now, according to~\cite[Ch.~III, 3.31 (1), p.~49]{Wat} we have
\begin{equation*}
   |J_\nu(x)| \leq \frac{2^{-\nu} x^\nu}{\Gamma(\nu+1)},
   \qquad \nu > -1/2.
\end{equation*}
Therefore,
\begin{align*}
   |\jna(x)|
   &= \sqrt{2(\alpha+2n+1)}\,|J_{\alpha+2n+1}(x)|
   \,x^{-\alpha-1} \\
   &\leq \frac{\sqrt{2(\alpha+2n+1)}\,2^{-(\alpha+2n+1)} x^{2n}}
             {\Gamma(\alpha+2n+2)},
\end{align*}
so that
\begin{equation}
\label{eq:summand}
   |\cna(f) \jna(x)| \leq C \lpnorm{f}p{\mua}
   \,\frac{n^{\delta+1/2} (x/2)^{2n}}{\Gamma(\alpha+2n+2)}
\end{equation}
and the series $\sum_{n=0}^\infty \cna(f) \jna(x)$ converges absolutely.
\end{proof}

\begin{lem}
\label{lem:deriv}
Let $\alpha > -1$, $p_0(\alpha)<p< p_1(\alpha)$,
and $\{\mu_n\}_{n=0}^\infty$ a sequence of positive numbers such that, for some positive constant $c$, $\mu_n \ge cn$ for every~$n$. Then, for any $f \in L^p(\mua)$ the series $\sum_{n=0}^\infty r^{\mu_n} \cna(f)\jna(x)$ is infinitely differentiable with respect to both variables $x \in (0,\infty)$ and~$r \in (0,1)$.
\end{lem}

\begin{proof}
For every $z \in (0,+\infty)$ and every $s \in (0,1)$, we can apply the arguments in the proof of Lemma~\ref{lem:convae} and get
\[
   | \mu_n r^{\mu_n-1} \cna(f)\jna(x) |
   \leq C \lpnorm{f}p{\mua} \mu_n s^{\mu_n-1}
   \,\frac{n^{\delta+1/2} (z/2)^{2n}}{\Gamma(\alpha+2n+2)}
\]
uniformly for $x \in (0,z)$, $r \in (0,s)$. The series
$
   \sum_{n=0}^\infty \mu_n s^{\mu_n-1}
   \,\frac{n^{\delta+1/2} (z/2)^{2n}}{\Gamma(\alpha+2n+2)}
$
is easily seen to be convergent, so that, by the dominated convergence theorem, the series $\sum_{n=0}^\infty r^{\mu_n} \cna(f)\jna(x)$ is differentiable with respect to $r$, its derivative is the term-by-term differentiated series
$
   \sum_{n=0}^\infty \mu_n r^{\mu_n-1} \cna(f)\jna(x)
$,
and this is a continuous function.

In order to prove the result for the first derivative with respect to $x$ we observe that, due to the definition of the functions $\jna$, it is enough to prove  that the series $\sum_{n=0}^\infty r^{\mu_n} \cna(f) \sqrt{2(2n+\alpha+1)}\, J_{2n+\alpha+1}(x)$ can be differentiated term by term. For this purpose we recall the formula $2J'_{\nu}(z) = J_{\nu-1}(z)-J_{\nu+1}(z)$, then again we can apply the arguments in the proof of the last lemma. 
The derivatives of higher order are handled in the same way.
\end{proof}

\section{Boundedness of the heat and Poisson semigroups}

Let us introduce some notation.
Given a sequence $\{a_n\}$, we will denote $\Delta a_n = a_n - a_{n-1}$.
Assume that the series $\sum_{n=0}^\infty a_n b_n$ and $\sum_{n=0}^\infty a_{n+1} b_n$ are convergent and $b_{-1}=0$, then it is easy to check that
\begin{equation}
\label{eq:anbn}
  \sum_{n=0}^\infty a_n \Delta b_n = - \sum_{n=0}^\infty \Delta a_{n+1} b_n.
\end{equation}

Given a function $f$ we shall denote by $\Sna f(x)$ its Fourier series 
$$
  \Sna f(x) = \sum_{k=0}^n \cka(f)\jka(x)
$$
where $\cka(f)$ are defined in~\eqref{eq:fourier}.
We shall also consider the Ces\`aro means of order one, defined as
\begin{equation}
\label{eq:Cesaro}
  C_n^\alpha f = \frac{S_0^\alpha f + S_1^\alpha f + \cdots + S_n^\alpha f}{n+1}.
\end{equation}

We shall use the following result whose proof can be found in \cite{CiStV}:

\begin{theor}
\label{th:Cesaro}
Let $\alpha > -1$ and $p_0(\alpha) < p < p_1(\alpha)$.
Then,
\begin{equation*}
  \lpnorm{C_n^\alpha f}{p}{\mua} \le C \lpnorm{f}{p}{\mua}
\end{equation*}
with a constant $C$ independent of~$n$.
\end{theor}

\begin{rem}
Actually, the Ces\`aro means are not directly studied in~\cite{CiStV}.
Instead, a different summation method
$$
  R_n^{\alpha}f = \frac{\rho_0 S_0^{\alpha}f + \cdots
  + \rho_n \Sna f}{\rho_0 + \cdots + \rho_n}
$$
(with $\rho_k = 2(\alpha+2k+2)$) is used.
But, as established in that paper, this method is equivalent to the given by the Ces\`aro means of order one, so the uniform boundedness of $R_n^\alpha$ is equivalent to the uniform boundedness of~$C_n^\alpha$.
\end{rem}

\begin{propo}
\label{pro:Cesaro}
Let $\alpha > -1$ and $\{\mu_n\}_{n=0}^\infty$ be a sequence of positive numbers such that, for some positive constant $c$, $\mu_n \ge cn$ for every~$n$. Given a function $f \in L^p(\mua)$ with $p_0(\alpha) < p < p_1(\alpha)$,
we shall consider the series 
$$
  V_r f(x) = \sum_{n=0}^\infty r^{\mu_n} \cna(f) \jna(x),
  \quad 0 < r < 1.
$$
Then for each $r$, $0< r < 1$, the series is absolutely convergent and we have 
$$
  V_r f(x) = \sum_{n=0}^\infty (\Delta^2 r^{\mu_{n+2}}) (n+1) C_n^\alpha f(x).
$$
\end{propo}

\begin{proof}
By Lemma~\ref{lem:convae}, $\sum_{n=0}^\infty |\cna(f) \jna(x)|$ converges for every $x \in (0,\infty)$, so their partial sums are bounded. Then, it is clear that there exists some positive number $t(x)$ such that $|\Sna f(x)| \le t(x)$ $\forall n$. As a consequence,
$$
  \sum_{n=0}^\infty r^{\mu_n} |\Sna f(x)|
  \le t(x) \sum_{n=0}^\infty r^{cn} < \infty
$$
for every $x \in (0,\infty)$ and every $r \in (0,1)$, and the same happens
with $\sum_{n=0}^\infty r^{\mu_{n+1}} |\Sna f(x)|$.
Thus, we can apply~\eqref{eq:anbn} and get 
\begin{align*}
  V_r f(x) &= \sum_{n=0}^\infty r^{\mu_n} \cna(f) \jna(x)
  = \sum_{n=0}^\infty r^{\mu_n} (S_n^\alpha f(x) - S_{n-1}^\alpha f(x)) \\
  &= - \sum_{n=0}^\infty \Delta r^{\mu_{n+1}} \Sna f(x).
\end{align*}
Let us show now that 
\begin{equation}
\label{eq:twoCn}
  \sum_{n=0}^\infty \Delta r^{\mu_{n+1}} (n+1) C_n^\alpha f(x)
\quad
\text{and}
\quad
  \sum_{n=0}^\infty \Delta r^{\mu_{n+1}} n C_{n-1}^\alpha f(x)
\end{equation}
are convergent series. We have already seen that $|\Sna f(x)| \le t(x)$ $\forall n$. Consequently, also $|C_n^\alpha f(x)| \le t(x)$ $\forall n$. Since $\mu_n \ge cn$, we have
$$
  |\Delta r^{\mu_{n+1}}| \le r^{\mu_{n+1}} + r^{\mu_n}
  \le r^{c(n+1)} + r^{cn} = (r^c + 1) r^{cn}.
$$
With this, 
\[
  \sum_{n=0}^\infty |\Delta r^{\mu_{n+1}}| (n+1) |C_n^\alpha f(x)|
  \le t(x) (1+r^c) \sum_{n=0}^\infty r^{cn} (n+1),
\]
which is convergent for every $x \in (0,\infty)$ and every $r \in (0,1)$.
The second series in~\eqref{eq:twoCn} can be analyzed analogously.
Hence as $\Sna f = (n+1) C_n^\alpha f - n C_{n-1}^\alpha f$,
\eqref{eq:anbn} can be used and we have
\begin{align*}
  V_r f(x) &= - \sum_{n=0}^\infty \Delta r^{\mu_{n+1}}
    \left( (n+1) C_n^\alpha f(x) - n C_{n-1}^\alpha f(x) \right) \\
  &= - \sum_{n=0}^\infty \Delta r^{\mu_{n+1}}
    \Delta\left( (n+1) C_n^\alpha f(x)\right) \\
  &= \sum_{n=0}^\infty (\Delta^2 r^{\mu_{n+2}}) (n+1) C_n^\alpha f(x).
  \qedhere
\end{align*}
\end{proof}

\begin{theor}
\label{th:Poisson}
Let $\alpha > -1$ and $p_0(\alpha) < p < p_1(\alpha)$.
For each function $f\in L^p(\mua)$ and for each $r$, $0<r<1$,
the function
$$
  P_r f(x) = \sum_{n=0}^\infty r^{\alpha+2n+1} \cna(f) \jna(x)
$$
is well defined. 
Moreover there exists a constant $C$, independent of $f$ and $r$, such that
$$
  \lpnorm{P_r f}{p}{\mua} \le C r^{\alpha+1} \lpnorm{f}{p}{\mua}.
$$
\end{theor}

\begin{proof}
We apply Proposition~\ref{pro:Cesaro} with 
$\mu_n = \alpha+2n+1$, then $\Delta^2 r^{\mu_{n+2}} = r^{\mu_n} (r^2 - 1)^2 \ge 0$.
Hence by using Theorem~\ref{th:Cesaro} we get
  \begin{align*}
  \lpnorm{P_r f}{p}{\mua}
  &= \Biglpnorm{\sum_{n=0}^\infty (\Delta^2 r^{\mu_{n+2}}) (n+1) C_n^\alpha f}{p}{\mua} \\
  &\le \sum_{n=0}^\infty (\Delta^2 r^{\mu_{n+2}}) (n+1) \lpnorm{C_n^\alpha f}{p}{\mua} \\
  &\le C \sum_{n=0}^\infty (\Delta^2 r^{\mu_{n+2}}) (n+1) \lpnorm{f}{p}{\mua} \\
  &= C r^{\mu_0} \lpnorm{f}{p}{\mua},
\end{align*}
because
\[
  \sum_{n=0}^\infty (\Delta^2 r^{\mu_{n+2}}) (n+1)
  = -\sum_{n=0}^\infty \Delta r^{\mu_{n+1}} \Delta (n+1)
  = -\sum_{n=0}^\infty \Delta r^{\mu_{n+1}} = r^{\mu_0}.
  \qedhere
\]
\end{proof}

In the case of the heat semigroup $W_r f$ we must apply more delicate arguments. The reason is that for any fixed $r$ the coefficients $\Delta^2 r^{\mu_{n+2}}$ take both signs. 

\begin{theor}
\label{th:heat}
Let $\alpha > -1$ and $p_0(\alpha) < p < p_1(\alpha)$.
For each function $f\in L^p(\mua)$ and for each $r$, $0<r<1$, the function
$$
  W_r f(x) = \sum_{n=0}^\infty r^{(\alpha+2n+1)^2} \cna(f) \jna(x)
$$
is well defined. 
Moreover there exists a constant $C$, independent of $f$ and $r$, such that
$$
  \lpnorm{W_r f}{p}{\mua} \le C \lpnorm{f}{p}{\mua}.
$$
\end{theor}

\begin{proof}
We can assume that $1/2 < r < 1$, for the case $0 < r \leq 1/2$ can be easily handled with the arguments of Lemma~\ref{lem:convae}. Now, it is easy to see that
\begin{align*}
  \sum_{n=N}^\infty (\Delta^2 r^{\mu_{n+2}}) (n+1)
  &= r^{\mu_{N+1}} - (N+1) \Delta r^{\mu_{N+1}},
  \\
  \sum_{n=0}^{N-1} (\Delta^2 r^{\mu_{n+2}}) (n+1)
  &= r^{\mu_0} - r^{\mu_{N+1}} + (N+1) \Delta r^{\mu_{N+1}}
\end{align*}
for any positive integer~$N$. Let us investigate the sign of $\Delta^2 r^{\mu_{n+2}}$. We have $\mu_n = (\alpha+2n+1)^2$, so that
\[
   \Delta^2 r^{\mu_{n+2}} = r^{\mu_n}
   (r^{\mu_{n+2}-\mu_n} - 2 r^{\mu_{n+1}-\mu_n} + 1)
   = r^{\mu_n} (r^8 s^2 - 2 s + 1),
\]
with $s = r^{4(\alpha + 2 n + 2)}$. Therefore
\[
   \Delta^2 r^{\mu_{n+2}} = 0
   \iff
   s = \frac{1}{1 + \sqrt{1 - r^8}}
\]
(the other solution does not belong to the interval $(0,1)$), that is,
\[
   4(\alpha + 2 n + 2) = \frac{\log (1 + \sqrt{1 - r^8})}{-\log r} 
   \sim (1-r)^{-1/2}.
\]
Here, $a(r) \sim b(r)$ means $C_1 \leq a(r) / b(r) \leq C_2$ for some positive constants $C_1$, $C_2$ independent of $r \in (1/2,1)$. This proves that there exists some $N(r) \sim (1-r)^{-1/2}$ such that
\[
  \Delta^2 r^{\mu_n+2}
  \begin{cases}
    < 0, & \text{when } n < N(r), \\
    \ge 0, & \text{when } n \ge N(r).
  \end{cases}
\]
By Proposition~\ref{pro:Cesaro}, the series that defines $W_r f(x)$ is absolutely convergent and moreover we have
\begin{align*}
  \lpnorm{W_r f}{p}{\mua}
  &= \Biglpnorm{\sum_{n=0}^\infty (\Delta^2 r^{\mu_{n+2}}) (n+1)
    C_n^\alpha f}{p}{\mua} \\
  &\le \sum_{n=0}^{N(r)-1} (-\Delta^2 r^{\mu_{n+2}}) (n+1)
    \lpnorm{C_n^\alpha f}{p}{\mua} \\
  &\qquad\qquad  + \sum_{n=N(r)}^\infty (\Delta^2 r^{\mu_{n+2}}) (n+1)
    \lpnorm{C_n^\alpha f}{p}{\mua}\\
  &\le C \sum_{n=0}^{N(r)-1} (-\Delta^2 r^{\mu_{n+2}})
    (n+1) \lpnorm{f}{p}{\mua} \\
  &\qquad\qquad + C \sum_{n=N(r)}^\infty (\Delta^2 r^{\mu_{n+2}}) (n+1)
    \lpnorm{f}{p}{\mua} \\
  &= C (-r^{\mu_0} + 2 r^{\mu_{N(r)+1}} - 2 (N(r)+1) \Delta r^{\mu_{N(r)+1}})
    \lpnorm{f}{p}{\mua}.
\end{align*}
Finally, the estimate $N(r) \sim (1-r)^{-1/2}$ gives
\begin{align*}
   \mu_{N(r)+1} - \mu_{N(r)} = 4 (\alpha + 2 N(r) + 2) &\sim (1-r)^{-1/2},
   \\
   (\log r) (\mu_{N(r)} - \mu_{N(r)+1}) &\sim (1-r)^{1/2},
   \\
   - \Delta r^{\mu_{N(r)+1}} 
   = r^{\mu_{N(r)+1}} (r^{\mu_{N(r)} - \mu_{N(r)+1}} - 1)
   &\sim (1-r)^{1/2},
\end{align*}
and
\[
   -r^{\mu_0} + 2 r^{\mu_{N(r)+1}} - 2 (N(r)+1) \Delta r^{\mu_{N(r)+1}}
   \le C.
   \qedhere
\]
\end{proof}

\begin{rem}
The definition of the Poisson and heat semigroup as a series, like
in Theorems~\ref{th:Poisson} and~\ref{th:heat}, is not always possible when
other orthogonal systems are used. For instance, in a analogous study
for Hermite and Laguerre expansions, the corresponding series can be divergent,
although the Poisson and heat semigroups can be defined with an appropriate kernel; see~\cite[Lemmas~2 and~4]{Mu}.
\end{rem}

Note that Theorem~\ref{th:heat} and the subordination formula give
$$
   \| P_r f \|_{L^p(d\mu_\alpha)} \le C\| f \|_{L^p(d\mu_\alpha)},
$$
if $\alpha>-1$, $p_0(\alpha) < p < p_1(\alpha)$, $f\in L^p(d\mu_\alpha)$ and
$0<r<1$. This result is weaker than Theorem 2 and it is insufficient to obtain the boundedness of the fractional integral as in Theorem~\ref{th:frac} below.

\section{Mean and almost everywhere convergence}

\begin{theor}
\label{th:convLpVr}
Let $\alpha > -1$ and $p_0(\alpha) < p< p_1(\alpha)$. Then
$P_r f \to f$ and $W_r f \to f$ in the $L^p(\mua)$-norm when
$r \to 1^-$ for every $f \in B_{p,\alpha}$.
\end{theor}

\begin{proof}
Let $V_rf$ denote either $P_r f$ or $W_r f$. Let us recall that, as
$p_0(\alpha) < p< p_1(\alpha)$, then $\jna \in L^p(\mua)$, for every $n = 0, 1,2, \dots$. Then, it is clear that
$$
  \lpnorm{V_r\jna - \jna}{p}{\mua} = (1-r^{\mu_n})\lpnorm{\jna}{p}{\mua} \to 0
$$
when $r\to 1^-$. From this, it follows that 
$$
  \lpnorm{V_rg-g}{p}{\mua} \to 0,
$$
when $r\to 1^-$, for every $g\in \spa\{\jna\}_{n=0}^\infty$.

Now, let $f\in B_{p,\alpha}$ and $\varepsilon>0$.  We choose $g\in \spa\{\jna\}_{n=0}^\infty$ such that $\lpnorm{f-g}{p}{\mua} < \varepsilon$.
Then, we have
\begin{multline*}
  \lpnorm{V_rf-f}{p}{\mua} \\
  \le \lpnorm{V_r(f-g)}{p}{\mua} + \lpnorm{V_rg-g}{p}{\mua}
    + \lpnorm{f-g}{p}{\mua} \\
  \le (C+1)\varepsilon + \lpnorm{V_rg-g}{p}{\mua}.
\end{multline*}
As $\lpnorm{V_rg-g}{p}{\mua} \to 0$ when $r\to 1^-$, we conclude the proof.
\end{proof}

It is known, see \cite{V}, that in the case 
$\alpha \ge -1/2$, 
$p_0(\alpha) < p < p_1(\alpha)$,
the space $B_{p,\alpha}$ coincides with the space
\begin{equation}
\label{eq:Epa}
  \{ f \in L^p(\mua): \ma f = f\} = \ma(L^p(\mua)),
\end{equation}
where $\ma$ is the multiplier of $[0,1]$ for
the so-called modified Hankel transform of order~$\alpha$.

As a consequence, we have the following result:

\begin{theor}
\label{th:convMaf}
Let $\alpha \ge -1/2$ and $p_0(\alpha) < p < p_1(\alpha)$.
Then, $P_r f \to \ma f$ and $W_r f \to \ma f$ in the $L^p(\mua)$-norm when $r \to 1^-$
for every $f \in L^p(\mua)$.
\end{theor}

\begin{proof}
The multiplier $\ma$ satisfies $\ma(\ma f) = \ma(f)$,
so $\ma f$ belongs to the set defined in~\eqref{eq:Epa} for every $f \in L^p(\mua)$
(see~\cite{V} for details).
Another important property of $\ma$ is
$$
  \int_0^\infty \ma f(y) g(y) \mua(y)
  = \int_0^\infty \ma g(y) f(y) \mua(y)
$$
for $f \in L^p(\mua)$ and $g \in L^{p'}(\mua)$.
In particular,
\begin{align*}
  \cna(\ma f)
  &= \int_0^\infty (\ma f) \jna \,\mua
  = \int_0^\infty (\ma \jna) f \,\mua \\
  &= \int_0^\infty f \jna \,\mua = \cna(f).
\end{align*}
Then, $P_r f = P_r (\ma f)$ and $W_r f = W_r (\ma f)$,
so the proof follows by applying Theorem~\ref{th:convLpVr}
to $\ma f$.
\end{proof}


Now we shall deal with the pointwise convergence of the heat and Poisson semigroups.
First we state a proposition which is parallel to Proposition~\ref{pro:Cesaro}.

\begin{propo}
\label{pro:convaeVr}
Let $\alpha > -1$, $p_0(\alpha) < p < p_1(\alpha)$ and $\{\mu_n\}_{n=0}^\infty$
be a sequence of positive numbers such that, for some positive constant $c$,
$\mu_n \ge cn$ for every~$n$. Then, for any $f \in L^p(\mua)$, 
$$
  \lim_{r \to 1^-} \sum_{n=0}^\infty r^{\mu_n} \cna(f)\jna(x)
  = \sum_{n=0}^\infty \cna(f)\jna(x)
$$
for every $x \in (0,\infty)$.
\end{propo}

\begin{proof}
We make a slight modification of the proof of Lemma~\ref{lem:convae}.
For a given $f$ and every $x \in (0,\infty)$ fixed, let us take
$g_r(n) = r^{\mu_n} \cna(f)\jna(x)$.
With this notation, \eqref{eq:summand} shows that
\begin{equation}
  |g_r(n)| \leq C r^{\mu_n} \lpnorm{f}p{\mua}
  \,\frac{n^{\delta+1/2} (x/2)^{2n}}{\Gamma(\alpha+2n+2)}
\end{equation}
so, taking
\begin{equation}
  g(n) = C \lpnorm{f}p{\mua}
  \,\frac{n^{\delta+1/2} (x/2)^{2n}}{\Gamma(\alpha+2n+2)}
\end{equation}
we have $|g_r(n)| \le g(n)$ for every $r \in (0,1)$,
and $\sum_{n=0}^\infty g(n) < \infty$.
Then, the dominated convergence theorem gives the result.
\end{proof}

Now, we can already state the following result:

\begin{theor}
\label{th:convaeVr}
Let $\alpha > -1$ and $p_0(\alpha) < p < p_1(\alpha)$.
Then, $P_r f$ and $W_r f$ converge almost everywhere to $f$
when $r \to 1^-$ for every $f \in B_{p,\alpha}$.
\end{theor}

\begin{proof}
Let $V_r$ denote either $P_r$ or $W_r$. Theorem~\ref{th:convLpVr} shows that under these conditions, $V_r f \to f$ in the $L^p(\mua)$-norm. Consequently, there exists a subsequence $\{r_j\}_{j=0}^\infty$ such that $V_{r_j} f \to f$ almost everywhere. This, in conjunction with Proposition~\ref{pro:convaeVr}, proves the theorem.
\end{proof}

Finally, in a similar way to Theorem~\ref{th:convMaf}, we have

\begin{theor}
\label{th:convaeMaf}
Let $\alpha \ge -1/2$ and $p_0(\alpha) < p < p_1(\alpha)$.
Then, $\lim_{r \to 1^-} P_r f = \ma f$ and $\lim_{r \to 1^-} W_r f = \ma f$
almost everywhere for every $f \in L^p(\mua)$.
\end{theor}

\section{Applications}

\subsection{Heat and Poisson equations}

\begin{theor}
\label{th:DE}
Let $\alpha > -1$, $p_0(\alpha) < p < p_1(\alpha)$ and  $f$ be a function in $L^p(\mua)$.
Then the functions $w(x,t) = W_{e^{-t}}f(x)$ (see Theorem~\ref{th:heat})
and $u(x,t) = P_{e^{-t}}f(x)$ (see Theorem~\ref{th:Poisson}) are infinitely differentiable in both variables $x\in (0,\infty)$ and $t\in(0,\infty)$.
They satisfy the differential equations
$\left( \frac{\partial}{\partial t} + L_\alpha \right) w(x,t) = 0$
and $\left( \frac{\partial^2}{\partial^2 t} + L_\alpha \right) u(x,t) = 0$.
Moreover, if $f\in B_{p,\alpha}$, the functions $w$ and $u$ are, respectively,
the solutions of the initial value problems given by the above differential equations
with the initial condition $w(x,0)=f(x)$ and $u(x,0)=f(x)$.
\end{theor}

\begin{proof}
As a direct consequence of Lemma~\ref{lem:deriv} we get the differentiability
of the functions $w(x,t)$ and $u(x,t)$.
Hence they satisfy the  corresponding differential equations.
In order to finish the proof we use Theorem~\ref{th:convaeVr}.
\end{proof}

\subsection{Fractional integral and potential spaces}

Let us consider the fractional integral of order $\lambda > 0$, $L_\alpha^{-\lambda/2}$, given by~\eqref{eq:frac} and~\eqref{eq:Poisson}. This
operator can be defined also by 
\[
   L_\alpha^{-\lambda/2} f
   =\sum_{k=0}^\infty c_k^\alpha (\alpha+2n+1)^{-\lambda} j_n^\alpha.
\]
By using Lemma 1 we can obtain the boundedness of the operator
$L_\alpha^{-\lambda/2}$ for some values of $\lambda$ depending on $\alpha$
and~$p$. However, as we show now, Theorem~\ref{th:Poisson} and the representation formula~\eqref{eq:Poisson} allow us to prove the boundedness of the fractional integral $L_\alpha^{-\lambda/2}$ on $L^p(d\mu_\alpha)$ for every $\lambda>0$.

After the change of variable $r=e^{-t}$, formula~\eqref{eq:Poisson} becomes
\begin{equation*}
  L_\alpha^{-\lambda/2}f
  = \frac{1}{\Gamma(\lambda)} 
  \int_0^1 ({-\log r})^{\lambda-1} P_r f \,\frac{dr}{r}
\end{equation*}
(remember that, as usual, we still write $P_r$ instead of $P_{-\log r}$).

\begin{theor}
\label{th:frac}
Let $\alpha > -1$ and $p_0(\alpha) < p < p_1(\alpha)$.
Then, for every $\lambda > 0$, we get 
\begin{equation*}
  \lpnorm{L_\alpha^{-\lambda/2} f}{p}{\mua} \le C \lpnorm{f}{p}{\mua}
\end{equation*}
with a constant $C=C(p,\alpha,\lambda)$ independent of~$f$.
\end{theor}

\begin{proof}
The boundedness of the Poisson semigroup (Theorem~\ref{th:Poisson}) gives
\begin{align*}
  \lpnorm{L^{-\lambda/2} f(x)}{p}{\mua(x)}
  &\le \frac{1}{\Gamma(\lambda)}
    \int_0^1 \lpnorm{P_rf(x)}{p}{\mua(x)} ({-\log r})^{\lambda-1} \,\frac{dr}{r} \\
  &\le C \lpnorm{f(x)}{p}{\mua(x)}
    \int_0^1 r^{\alpha+1} ({-\log r})^{\lambda-1} \,\frac{dr}{r} \\
  &\le C' \lpnorm{f(x)}{p}{\mua(x)}.
  \qedhere
\end{align*}
\end{proof}

The above theorem has the following consequence. For $\alpha > -1$ and $s > 0$ define, as usual, the potential space $\calLpsa$ by (see \cite{Stein-2})
$$
  \calLpsa = \{\, f:  \exists g \in L^p(\mua)
    \text{ such that }  L_\alpha^{-s/2}g = f \,\}
  = L_\alpha^{-s/2}(L^p(\mua))
$$
with the norm $\calnorm{f} = \lpnorm{g}{p}{\mua}$.
Then, for $p_0(\alpha) < p < p_1(\alpha)$,
$$
  \lpnorm{f}{p}{\mua} = \lpnorm{L_\alpha^{-s/2}g}{p}{\mua}
  \le C \lpnorm{g}{p}{\mua} = C\calnorm{f}
$$
and so $\calLpsa \subset L^p(\mua)$.


\end{document}